\documentclass{paper}
\usepackage[utf8]{inputenc}
\usepackage{float}
\usepackage{amsmath}
\usepackage{amsthm}
\usepackage{amssymb}
\usepackage{graphicx}
\usepackage[numbers]{natbib}
\usepackage[unicode=true,
 bookmarks=false,
 breaklinks=false,pdfborder={0 0 1},backref=section,colorlinks=false]
 {hyperref}

\makeatletter

\newtheorem{cor}{Corollary}
\newtheorem{defn}{Definition}
\newtheorem{prop}{Definition}

\usepackage{stmaryrd}

\makeatother

\begin{document}
\title{Asymptotically minimal contractors based on the centered form;\\
Application to the stability analysis of linear systems}
\author{Luc Jaulin \thanks{Luc Jaulin is with the team Robex of Lab-STICC, ENSTA-Bretagne, Brest,
France, e-mail: \protect\href{http://lucjaulin@gmail.com}{lucjaulin@gmail.com}.}}
\maketitle
\begin{abstract}
This paper proposes a new interval-based contractor for nonlinear
equations which is minimal when dealing with narrow boxes. The method
is based on the centered form classically used by interval algorithms
combined with a Gauss Jordan band diagonalization preconditioning.
As an illustration in stability analysis, we propose to compute the
set of all parameters of a characteristic function of a linear dynamical
system which have at least one zero in the imaginary axis. Our approach
is able compute a guaranteed and accurate enclosure of the solution
set faster than existing approaches. 
\end{abstract}


\section{Introduction}

Interval analysis is an efficient tool used for solving rigorously
complex nonlinear problems involving bounded uncertainties \cite{Ceberio01}
\cite{Kreinovich:97} \cite{Rauh11}. Many interval algorithms are
based on the notion of \emph{contractor} \cite{ChabertJaulin09} which
is an operator which shrinks an axis-aligned box $[\mathbf{x}]$ of
$\mathbb{R}^{n}$ without removing any point of the solution set $\mathbb{X}$.
The set $\mathbb{X}$ is assumed to be defined by equations involving
the components $x_{1},\dots,x_{n}$ of a vector $\mathbf{x}\in\mathbb{R}^{n}$. 

Combined with a paver \cite{Raazesh09} which bisects boxes, the contractor
builds an outer approximation of the set $\mathbb{X}$. The resulting
methodology can be applied in several domains of engineering such
as identification \cite{RamdaniPoignet05}, localization \cite{Langerwisch:Wagner:2012}
\cite{Guyonneau2013}, SLAM \cite{mustafa18} \cite{robloc}, vision
\cite{DBLP:EhambramVW21}, reachability \cite{CdC20}, control \cite{Rauh2009IntervalAT}
\cite{JianWan07}, calibration \cite{Daney06}, etc.

\emph{Centered form} is one of the most fundamental brick in interval
analysis. It is traditionally used to enclose the range of a function
over narrow intervals \cite{Moore79}\cite{Neumaier90}\cite{Hansen92}.
The quadratic approximation property, guarantees an asymptotically
small overestimation for sufficiently narrow boxes. In this paper,
we propose to use the centered form to build efficient contractors
\cite{JaulinBook01} that are optimal when the intervals are narrow. 

To achieve this goal, we first get a guaranteed first order enclosure
of each equation composing our problem. Then, we combine these constraints
preserving the first order approximation using interval linear techniques.
More particularly, we propose to use a preconditioning method based
on a Gauss-Jordan band diagonalization. We show that our approach
is guaranteed to enclose all solutions of the problem and that it
outperforms state of the art techniques. 

The main contribution of this paper is that the contractor we propose
is asymptotically minimal, \emph{i.e.}, it is minimal when the boxes
are small. To the best of my knowledge, such a contractor does not
exist in the literature even if some use a linear approximation (see
the X-Taylor iteration \cite{Araya12} tested on global minimization
problems, \cite{BraemsScan00} which is similar to X-Taylor but for
solving inequalities, the interval Newton \cite{Moore79} used for
solving square nonlinear systems, or the affine arithmetic \cite{Figueiredo04}
which has been used for non-square systems but which is not asymptotically
minimal).

Section \ref{sec:preliminaries} recalls some useful mathematical
notions related to the sensitivity of the solution set of a linear
system. Section \ref{sec:wrappers} introduces wrappers to approximate
accurately a function over a box. Section \ref{sec:asymptotically-minimal-contractor}
defines what is an asymptotically minimal contractor and Section \ref{sec:centered-contractor}
gives an algorithm to generate it. The relevance and the efficiency
of our approach are shown in Section \ref{sec:testcase} on the stability
analysis of a linear differential equation with delays. Section \ref{sec:Conclusion}
concludes the paper.

\section{Preliminaries\label{sec:preliminaries}}

This section recalls some basic definitions and theorems related to
the sensitivity of the solution set of a linear system with respect
to small perturbations. They will be used later in the paper to define
the asymptotic minimality of our approximation for the solution set. 

\subsection{Proximity}

\begin{figure}[h]
\begin{centering}
\includegraphics[width=7cm]{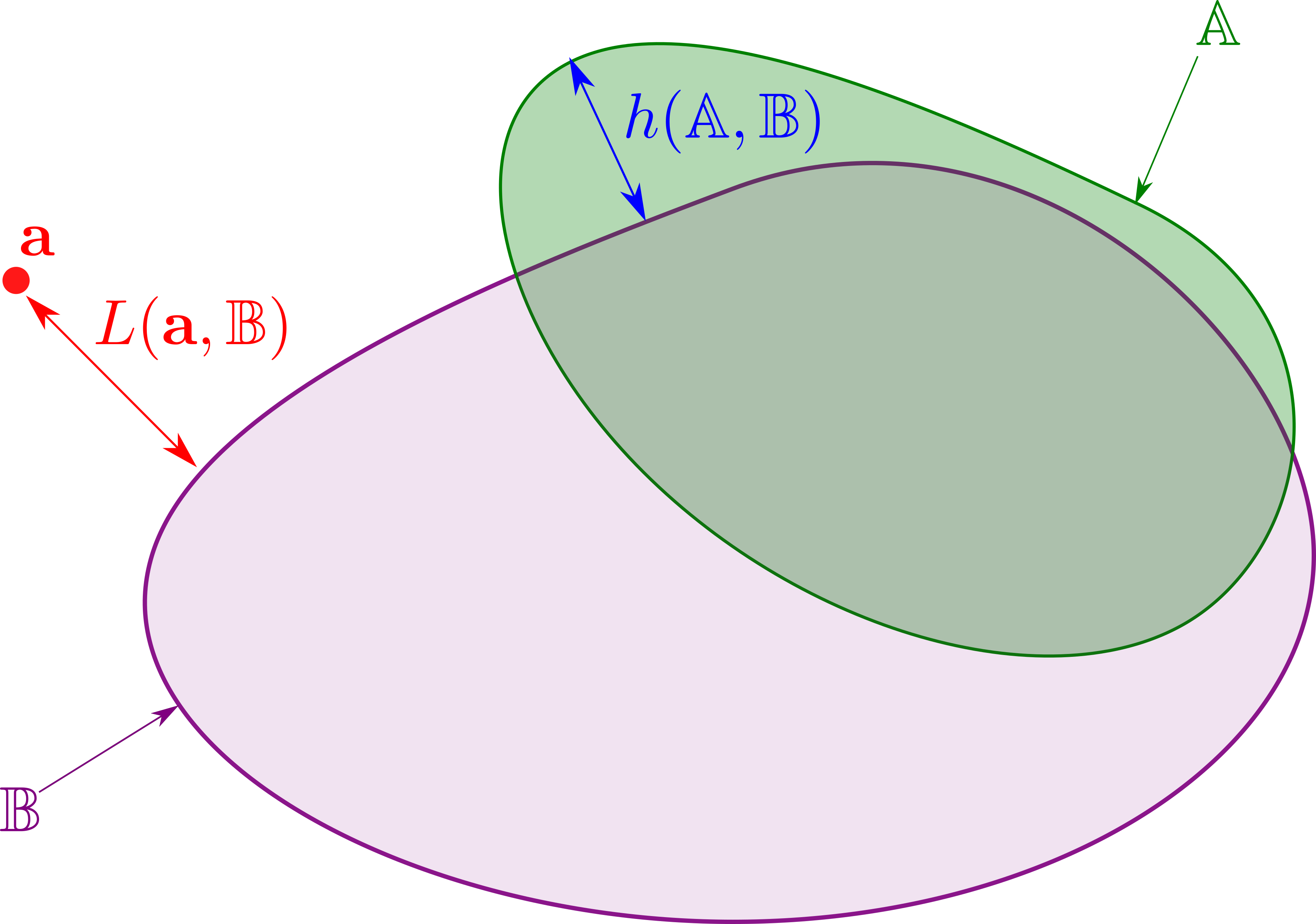}
\par\end{centering}
\caption{Proximity of $\mathbb{A}$ to $\mathbb{B}$}
\label{fig:norms.png}
\end{figure}

Denote by $L(\mathbf{a},\mathbf{b})$ the distance between $\mathbf{a}$
and $\mathbf{b}$ of $\mathbb{R}^{n}$ induced by the $L$-norm. As
illustrated by Figure \ref{fig:norms.png}, the \emph{proximity} of
$\mathbb{A}$ to $\mathbb{B}$, where $\mathbb{A}$ and $\mathbb{B}$
are closed subsets of $\mathbb{R}^{n}$, is defined by
\begin{equation}
h(\mathbb{A},\mathbb{B})=\sup_{\mathbf{a}\in\mathbb{A}}L(\mathbf{a},\mathbb{B})
\end{equation}
where
\begin{equation}
L(\mathbf{a},\mathbb{B})=\inf_{\mathbf{b}\in\mathbb{B}}L(\mathbf{a},\mathbf{b}).
\end{equation}
The norm $L$ that will be used later in the algorithm will be the
$L_{\infty}$ norm, even if, in the pictures, for a better visibility,
we use the Euclidean $L_{2}$ norm.

A nested sequence of closed subsets $\mathbb{B}(k)\subset\mathbb{R}^{n}$,
$k\in\mathbb{N}$ is converging to $\mathbf{x}$ if 
\begin{equation}
\begin{array}{c}
\lim_{k\rightarrow\infty}h(\mathbb{B}(k),\{\mathbf{x}\})=0.\end{array}
\end{equation}

\subsection{Linear systems}

The following proposition allows us to quantify the sensitivity of
the solutions of a linear system of equations.
\begin{prop}
\label{prop:sensitivity:lin}Consider a point $\mathbf{x}$ which
satisfies the linear system $\mathbf{A}\cdot\mathbf{x}=\mathbf{b}$,
where $\mathbf{A}$ has independent lines. Consider a small variation
$d\mathbf{A}$ of $\mathbf{A}$. The quantity 
\begin{equation}
d\mathbf{x}=-\mathbf{A}^{\dagger}\cdot(d\mathbf{A}\cdot\mathbf{x}+d\mathbf{A}\cdot d\mathbf{x})
\end{equation}
where
\begin{equation}
\mathbf{A}^{\dagger}=\mathbf{A}^{\text{T}}(\mathbf{A}\cdot\mathbf{A}^{\text{T}})^{-1}
\end{equation}
is the generalized inverse of $\mathbf{A}$, satisfies 
\begin{equation}
(\mathbf{A}+d\mathbf{A})\cdot(\mathbf{x}+d\mathbf{x})=\mathbf{b}.\label{eq:prop1}
\end{equation}
\end{prop}
This proposition tells us that if we move $\mathbf{A}$ a little,
then, the solution set for the linear equation moves a little also,
at order 1. 

\textbf{Proof}. We have 
\begin{equation}
\begin{array}{ccl}
 & \, & (\mathbf{A}+d\mathbf{A})\cdot(\mathbf{x}+d\mathbf{x})=\mathbf{b}\\
\Leftrightarrow &  & \mathbf{A}\cdot\mathbf{x}+\mathbf{A}\cdot d\mathbf{x}+d\mathbf{A}\cdot\mathbf{x}+d\mathbf{A}\cdot d\mathbf{x}=\mathbf{b}
\end{array}
\end{equation}
Thus 
\begin{equation}
\mathbf{A}\cdot d\mathbf{x}+d\mathbf{A}\cdot\mathbf{x}+d\mathbf{A}\cdot d\mathbf{x}=\mathbf{0}
\end{equation}
\emph{i.e.}
\begin{equation}
\mathbf{A}\cdot d\mathbf{x}=-d\mathbf{A}\cdot\mathbf{x}-d\mathbf{A}\cdot d\mathbf{x}
\end{equation}
Since $\mathbf{A}$ has independent lines, the solution which minimizes
$\|d\mathbf{x}\|$ is
\begin{equation}
d\mathbf{x}=\mathbf{A}^{\dagger}\cdot(-d\mathbf{A}\cdot\mathbf{x}-d\mathbf{A}\cdot d\mathbf{x}).\blacksquare
\end{equation}

\begin{cor}
\label{cor:moveplane}Consider the hyperplane 
\begin{equation}
\mathcal{P}=\{\mathbf{x}\in\mathbb{R}^{n}|\mathbf{A}\cdot\mathbf{x}=\mathbf{0}\},
\end{equation}
where $\mathbf{A}$ has independent lines. Consider a small variation
$d\mathbf{A}$ of $\mathbf{A}$ with $\|d\mathbf{A}\|=O(\varepsilon)$
where $\varepsilon$ is small. Take a point $d\mathbf{x}\in\mathcal{P}$
with $\|d\mathbf{x}\|=O(\varepsilon)$. The distance from $d\mathbf{x}$
to $\mathcal{\tilde{P}}=\{\mathbf{x}\in\mathbb{R}^{n}|(\mathbf{A}+d\mathbf{A})\cdot\mathbf{x}=\mathbf{0}\}$
is $o(\varepsilon)$, i.e., $O(\varepsilon^{2})$. 
\end{cor}
\textbf{Proof}. Denote by $\hat{\mathbf{p}}$ the projection of a
point $\mathbf{p}\in\mathcal{P}$ on $\mathcal{\tilde{P}}$. From
(\ref{eq:prop1}), we have 
\begin{equation}
\|\hat{\mathbf{p}}-\mathbf{p}\|=\|\mathbf{A}^{\dagger}\cdot d\mathbf{A}\cdot\mathbf{p}\|+o(\varepsilon).
\end{equation}
If we take $\mathbf{p}=d\mathbf{x}$. We get 
\begin{equation}
\begin{array}{ccl}
\|d\hat{\mathbf{x}}-d\mathbf{x}\| & = & \|\mathbf{A}^{\dagger}\cdot d\mathbf{A}\cdot d\mathbf{x}\|+o(\varepsilon)\\
 & = & o(\varepsilon)=O(\varepsilon^{2}).\blacksquare
\end{array}
\end{equation}

\section{Wrappers\label{sec:wrappers}}

\begin{figure}[h]
\begin{centering}
\includegraphics[width=7cm]{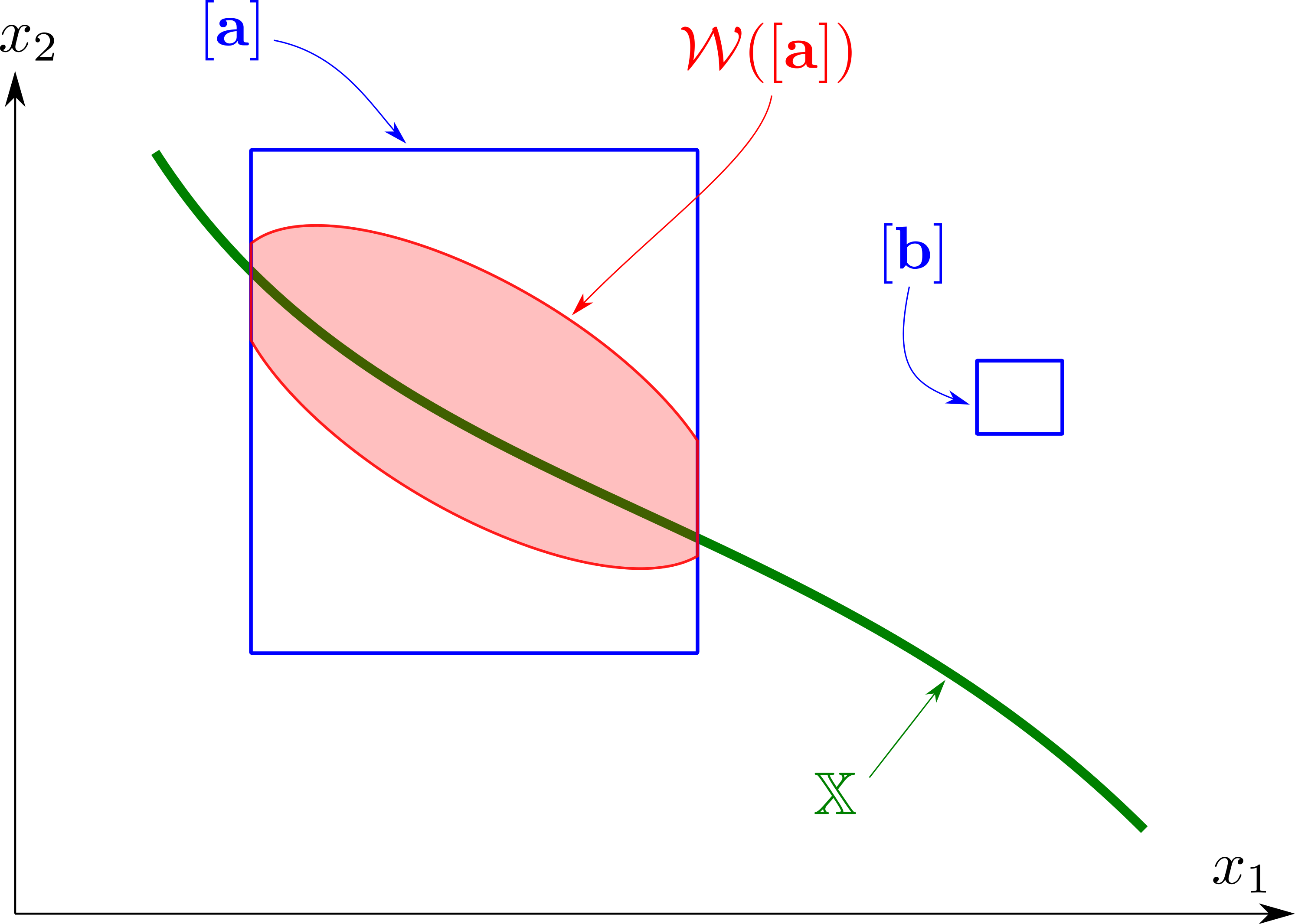}
\par\end{centering}
\caption{Illustration of a wrapper}
\label{fig:defwrapper}
\end{figure}

The approximation of sets using boxes computed using interval analysis
generates a strong wrapping effect. It has been shown by several authors
that it was possible to get a linear approximation with a better accuracy
using other types of sets such as zonotopes \cite{DBLP:Combastel15}
\cite{DBLP:Combastel22}, ellipsoids \cite{rauh:ellipse:2021}, or
doubleton \cite{DBLP:journals/cnsns/KapelaMWZ21}. Before defining
the notion of wrapper to quantify the order of approximation we can
get, we first recall what is a contractor.
\begin{defn}
Denote by $\mathbb{IR}^{n}$ the set of boxes of $\mathbb{R}^{n}$.
A \emph{contractor} associated to the closed set $\mathbb{X}\subset\mathbb{R}^{n}$
is a function $\mathcal{C}:\mathbb{IR}^{n}\mapsto\mathbb{IR}^{n}$
such that
\[
\begin{array}{ccc}
\mathcal{C}([\mathbf{x}])\subset[\mathbf{x}] & \, & \text{(contraction)}\\{}
[\mathbf{x}]\cap\mathbb{X}\subset\mathcal{C}([\mathbf{x}]) &  & \text{(consistency)}
\end{array}
\]
The contractor $\mathcal{C}$ for $\mathbb{X}$ is \emph{minimal}
if $\mathcal{C}([\mathbf{x}])=\llbracket[\mathbf{x}]\cap\mathbb{X}\rrbracket$
where $\llbracket\mathbb{A}\rrbracket$ denotes the smallest box enclosing
the set $\mathbb{A}$.
\end{defn}
The following definition of a \emph{wrapper }extends the concept of
contractor and will be needed for convergence analysis.
\begin{defn}
\label{def:wrappers}A \emph{wrapper} associated to the closed set
$\mathbb{X}\subset\mathbb{R}^{n}$ is a function $\mathcal{W}:\mathbb{IR}^{n}\mapsto\mathcal{P}(\mathbb{R}^{n})$
such that
\[
\begin{array}{ccc}
\mathcal{W}([\mathbf{x}])\subset[\mathbf{x}] & \, & \text{(contraction)}\\{}
[\mathbf{x}]\cap\mathbb{X}\subset\mathcal{W}([\mathbf{x}]) &  & \text{(consistency)}\\
\mathbf{x}\notin\mathbb{X}\Rightarrow\exists\varepsilon,\forall[\mathbf{x}]\subset B(\mathbf{x},\varepsilon),\mathcal{W}([\mathbf{x}])=\emptyset &  & \text{(accuracy)}
\end{array}
\]
where $B(\mathbf{x},\varepsilon)$ is the box with center $\mathbf{x}$
and radius $\varepsilon$.
\end{defn}
An illustration of a wrapper is given by Figure \ref{fig:defwrapper}.
The set $\mathbb{X}$ is a curve which could be given by an equation.
For the box $[\mathbf{a}]$, the set $\mathcal{W}([\mathbf{a}])$
encloses the part of $\mathbb{X}$ which is inside $[\mathbf{a}]$.
For the box $[\mathbf{b}]$, we have $\mathcal{W}([\mathbf{b}])=\emptyset$.

The wrapper $\mathcal{W}$ for $\mathbb{X}$ has an order $1$ at
point $\mathbf{x}$ if for all nested sequences of boxes $[\mathbf{x}](k)$
converging to $\mathbf{x}$, we have
\begin{equation}
\lim_{k\rightarrow\infty}\frac{h(\mathcal{W}([\mathbf{x}](k)),\mathbb{X})}{w([\mathbf{x}](k))}=0\label{eq:W:order1}
\end{equation}
where $w([\mathbf{x}])$ is the width of $[\mathbf{x}]$. Denote by
$\text{Wrap}(\mathbb{X},\mathbf{x})$ the set of all wrappers for
$\mathbb{X}$ which have an order $1$ at point $\mathbf{x}$.

\begin{figure}[h]
\begin{centering}
\includegraphics[width=7cm]{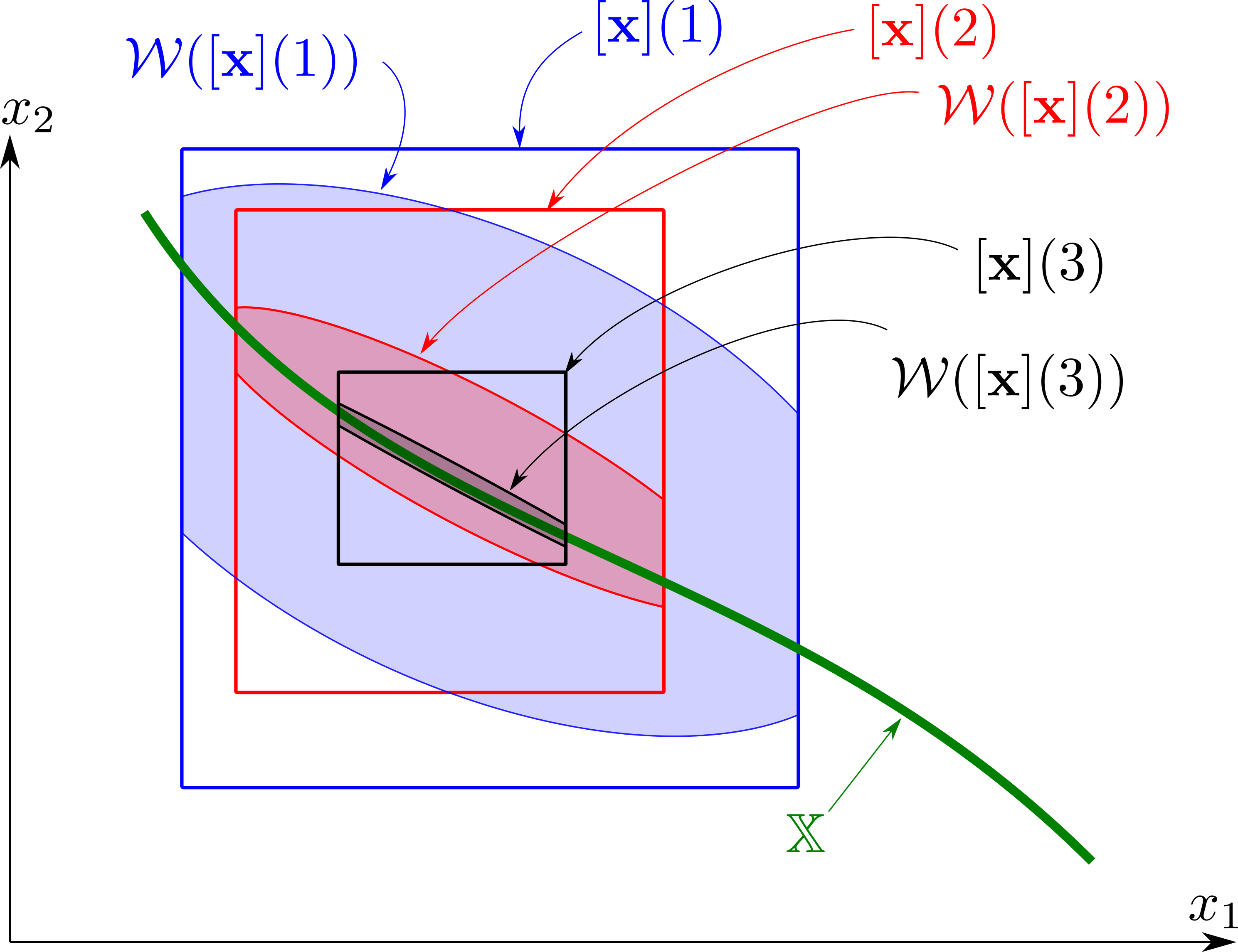}
\par\end{centering}
\caption{Wrapper of order 1}
\label{fig:converg_wrapper}
\end{figure}

The notion of order is illustrated by Figure \ref{fig:converg_wrapper}.
Larger is $k$, narrower is $[\mathbf{x}](k)$ and more accurate is
the approximation. 

\begin{figure}[h]
\begin{centering}
\includegraphics[width=7cm]{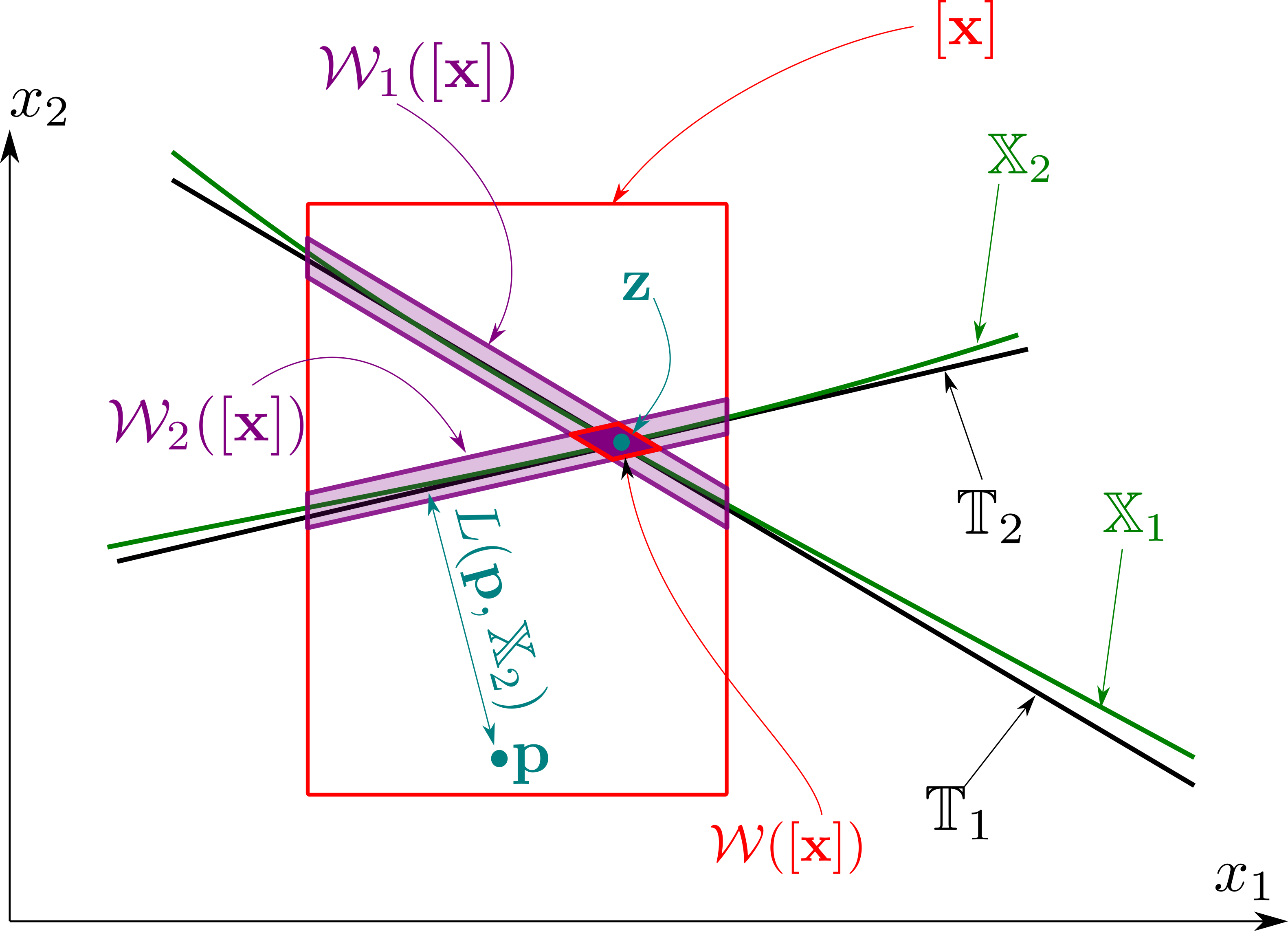}
\par\end{centering}
\caption{Intersection of two wrappers of order 1}
\label{fig:nwrapper}
\end{figure}

\begin{defn}
We define the intersection $\mathcal{W}$ of two wrappers $\mathcal{W}_{1}$
and $\mathcal{W}_{2}$ as
\begin{equation}
\mathcal{W}([\mathbf{x}])=(\mathcal{W}_{1}\cap\mathcal{W}_{2})([\mathbf{x}])=\mathcal{W}_{1}([\mathbf{x}])\cap\mathcal{W}_{2}([\mathbf{x}]).
\end{equation}
It is trivial to check that if $\mathcal{W}_{1}$ is a wrapper for
$\mathbb{X}_{1}$ and $\mathcal{W}_{2}$ is a wrapper for $\mathbb{X}_{2}$
then $\mathcal{W}=\mathcal{W}_{1}\cap\mathcal{W}_{2}$ is a wrapper
for $\mathbb{X}_{1}\cap\mathbb{X}_{2}.$ Unfortunately, the order
of the approximation is not always preserved. The following proposition
gives some conditions which allows us to preserve the order 1.
\end{defn}
\begin{prop}
\label{prop:n:wrappers}Given $m$ sets $\mathbb{X}_{i}=\left\{ \mathbf{x}\in\mathbb{R}^{n}|f_{i}(\mathbf{x})=0\right\} $,
where $f_{i}:\mathbb{R}^{n}\mapsto\mathbb{R}$. Consider $\mathbb{Z}=\bigcap_{i}\mathbb{X}_{i}$
and a point $\mathbf{z}\in\mathbb{Z}.$ Assume that all $\nabla f_{i}(\mathbf{z})$
are independent. If $\mathcal{W}=\bigcap_{i}\mathcal{W}_{i}$, we
have
\begin{equation}
\forall i,\mathcal{W}_{i}\in\text{Wrap}(\mathbb{X}_{i},\mathbf{z})\Rightarrow{\textstyle \bigcap}_{i}\mathcal{W}_{i}\in\text{Wrap}(\mathbb{Z},\mathbf{z})
\end{equation}
\end{prop}
Figure \ref{fig:nwrapper} illustrates that the intersection of two
wrappers of order 1 at $\mathbf{z}$ is generally a wrapper of order
1 at $\mathbf{z}$. In the figure, the set $\mathbb{Z}=\mathbb{X}_{1}\cap\mathbb{X}_{2}$
is the singleton $\{\mathbf{z}\}$. The box $[\mathbf{x}]$ should
be interpreted as a narrow box containing $\mathbf{z}$.

\textbf{Proof}. Since $\mathbb{Z}=\bigcap_{i}\mathbb{X}_{i}$, $\mathcal{W}=\bigcap_{i}\mathcal{W}_{i}$
is a wrapper for $\mathbb{Z}$. We also need to prove that the order
of $\mathcal{W}$ is 1 at $\mathbf{z}$. For this, consider a sequence
$[\mathbf{x}](k)$ converging to $\mathbf{z}$. When $k$ is large
$\varepsilon=w([\mathbf{x}](k))$ is small. For short, let us omit
the dependency with respect to $k$. For all $\mathbf{p}\in[\mathbf{x}],$
we have $\|\mathbf{p}-\mathbf{z}\|=O(\varepsilon)$. If $\mathbb{T}_{i}$
is the tangent space of $\mathbb{X}_{i}$ at point $\mathbf{z}$ then
\begin{equation}
L(\mathbf{p},\mathbb{X}_{i})=L(\mathbf{p},\mathbb{T}_{i})+o(\varepsilon).
\end{equation}
If all $\mathbb{T}_{i}$ are transverse, we have 
\begin{equation}
L(\mathbf{p},\mathbb{Z})=L(\mathbf{p},{\textstyle \bigcap}_{i}\mathbb{X}_{i})=L(\mathbf{p},{\textstyle \bigcap}_{i}\mathbb{T}_{i})+o(\varepsilon).\label{eq:loo_pZ}
\end{equation}
Take now, $\mathbf{p}\in\mathcal{W}([\mathbf{x}])$. Since $\forall i,L(\mathbf{p},\mathbb{T}_{i})=o(\varepsilon)$
and since the $\mathbb{T}_{i}$ are transverse, we get that $L(\mathbf{p},\bigcap_{i}\mathbb{T}_{i})=o(\varepsilon)$.
Therefore, from (\ref{eq:loo_pZ}), $L(\mathbf{p},\mathbb{Z})=o(\varepsilon)$.
Since this is true for all $\mathbf{p}\in\mathcal{W}([\mathbf{x}]),$
we have 
\begin{equation}
h(\mathcal{W}([\mathbf{x}]),\mathbb{Z})=\sup_{\mathbf{p}\in\mathcal{W}([\mathbf{x}])}L(\mathbf{p},\mathbb{Z})=o(\varepsilon)=o(w([\mathbf{x}])).
\end{equation}
Taking into account the dependency of $[\mathbf{x}]$ in $k$, we
get:
\begin{equation}
\lim_{k\rightarrow\infty}\frac{h(\mathcal{W}([\mathbf{x}](k)),\mathbb{Z})}{w([\mathbf{x}](k))}=0.\blacksquare
\end{equation}

\section{Asymptotically minimal contractor\label{sec:asymptotically-minimal-contractor} }

\begin{figure}[h]
\begin{centering}
\includegraphics[width=7cm]{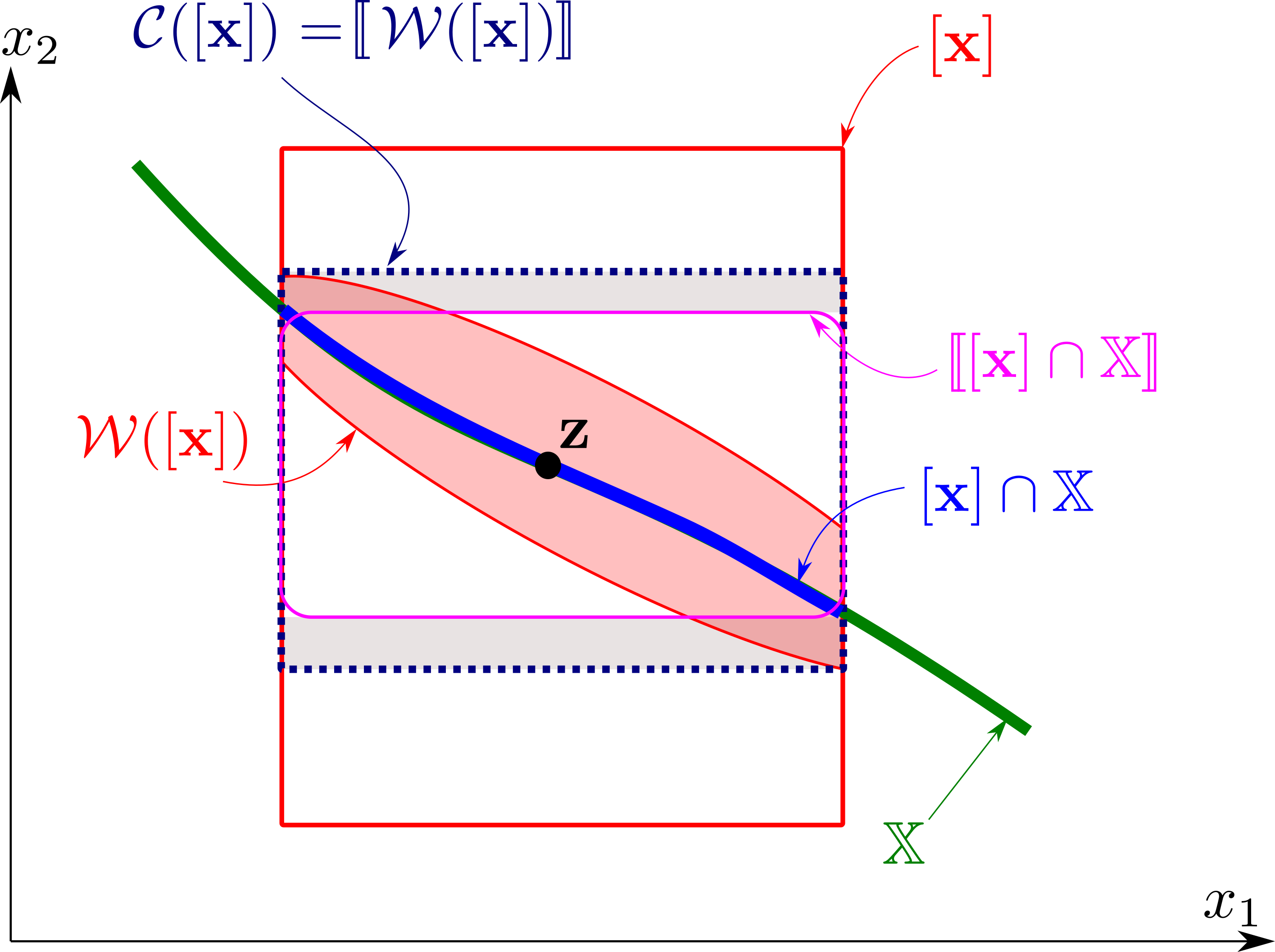}
\par\end{centering}
\caption{Asymptotic minimal contractor}
\label{fig:centered_asymp_ctc}
\end{figure}

Consider the special case where wrappers, as defined by Definition
\ref{def:wrappers}, generate sets $\mathcal{W}([\mathbf{x}])$ that
are boxes of $\mathbb{R}^{n}$. The order cannot be equal to 1 (it
can only be equal to $0$), except if $n=1$. Now, we can use the
wrappers of order 1, as an intermediate results, to get contractors
with a good accuracy. This section defines formally such accurate
contractors which is called \emph{asymptotically minimal}.
\begin{defn}
A contractor for $\mathbb{X}$ is \emph{asymptotically minimal} at
point $\mathbf{z}\in\mathbb{X}\subset\mathbb{R}^{n}$ if for any nested
sequence $[\mathbf{x}](k)$ converging to $\mathbf{z}$, we have 
\begin{equation}
\lim_{k\rightarrow\infty}\,\frac{h(\mathcal{C}([\mathbf{x}](k)),\llbracket[\mathbf{x}](k)\cap\mathbb{X}\rrbracket)}{w([\mathbf{x}](k))}=0.\label{eq:def:aminimal}
\end{equation}
\end{defn}
Note that since $\mathcal{C}$ is a contractor the quantity $\mathcal{C}([\mathbf{x}](k))$
is a box.
\begin{prop}
\label{prop:Cmini:fromW}If $\mathcal{W}\in\text{Wrap}(\mathbb{X},\mathbf{z})$,
then, the contractor defined by 
\begin{equation}
\mathcal{C}([\mathbf{x}])=\llbracket\mathcal{W}([\mathbf{x}])\rrbracket
\end{equation}
is an asymptotically minimal contractor for $\mathbb{X}$ at $\mathbf{z}$.
\end{prop}
An illustration of the proposition is given by Figure \ref{fig:centered_asymp_ctc}.
The gray part corresponds to the pessimism of the contractor which
tends to disappear when $[\mathbf{x}]$ becomes narrow.

\textbf{Proof}. The proof is by contradiction. Assume that $\mathcal{C}([\mathbf{x}])=\llbracket\mathcal{W}([\mathbf{x}])\rrbracket$
is not asymptotically minimal in $\mathbf{z}$. From (\ref{eq:def:aminimal}),
there exists a sequence of nested boxes such converging to $\mathbf{z}$
such that 
\begin{equation}
\lim_{k\rightarrow\infty}\,\frac{h(\llbracket\mathcal{W}([\mathbf{x}])(k)\rrbracket,\llbracket[\mathbf{x}](k)\cap\mathbb{X}\rrbracket)}{w([\mathbf{x}](k))}>0.
\end{equation}
Thus 
\begin{equation}
\lim_{k\rightarrow\infty}\,\frac{h(\mathcal{W}([\mathbf{x}])(k),\llbracket[\mathbf{x}](k)\cap\mathbb{X}\rrbracket)}{w([\mathbf{x}](k))}>0.
\end{equation}
Therefore
\begin{equation}
\lim_{k\rightarrow\infty}\frac{h(\mathcal{W}([\mathbf{x}](k)),\mathbb{X})}{w([\mathbf{x}](k))}>0.
\end{equation}
This is inconsistent with the fact that $\mathcal{W}$ has an order
1 in $\mathbf{z}$ (see (\ref{eq:def:aminimal})).$\blacksquare$

\section{Centered contractor \label{sec:centered-contractor}}

In this section, we show how to build an asymptotic minimal contractor
using the centered form. We will consider functions $\mathbf{f}:\mathbb{R}^{n}\mapsto\mathbb{R}^{p}$
which are all continuous and differentiable. More precisely, the functions
$\mathbf{f}$ are described by continuous operator of functions such
as $+,-,/,\sin,\exp,\dots$ As a consequence using interval analysis,
we are able to enclose the range of $\mathbf{f}$ and of $\frac{d\mathbf{f}}{d\mathbf{x}}$
over a box $[\mathbf{x}]$. In \cite{Moore79}, Moore has proved that
if $w([\mathbf{x}])=O(\varepsilon)$ then using interval computation,
we get an enclosure $[\mathbf{f}]([\mathbf{x}])$ for $\mathbf{f}([\mathbf{x}])$
and an enclosure $[\frac{d\mathbf{f}}{d\mathbf{x}}]([\mathbf{x}])$
for $\frac{d\mathbf{f}}{d\mathbf{x}}([\mathbf{x}])$ such that $w([\mathbf{f}]([\mathbf{x}]))=O(\varepsilon)$
and $w(\frac{d\mathbf{f}}{d\mathbf{x}}([\mathbf{x}]))=O(\varepsilon)$.

\subsection{Scalar case}
\begin{prop}
Consider the equation $f(\mathbf{x})=\mathbf{0}$, where \textbf{$f:\mathbb{R}^{n}\mapsto\mathbb{R}$}
is differentiable. The solution set is
\begin{equation}
\begin{array}{ccl}
\mathbb{X} & = & \{\mathbf{x}\in\mathbb{R}^{n}\,|\,f(\mathbf{x})=\mathbf{0}\}.\end{array}
\end{equation}
Consider a point $\mathbf{z}$ such that $f(\mathbf{z})=0$. Consider
a nested sequence $[\mathbf{x}](k)$ converging to $\mathbf{z}$.
The function $\mathcal{L}:\mathbb{IR}^{n}\mapsto\mathcal{P}(\mathbb{R}^{n})$
defined by
\begin{equation}
\begin{array}{ccl}
\mathcal{L}([\mathbf{x}]) & =\{ & \mathbf{x}\in[\mathbf{x}]\,|\,\exists\mathbf{a}\in[\frac{df}{d\mathbf{x}}]([\mathbf{x}]),\\
 &  & f(\mathbf{m})+\mathbf{a}\cdot(\mathbf{x}-\mathbf{m})=0,\\
 &  & \mathbf{m}=\text{center}([\mathbf{x}])\}
\end{array}\label{eq:Lxx}
\end{equation}
is a wrapper of order 1, i.e., it belongs to $\text{Wrap}(\mathbb{X},\mathbf{z})$.
It will be called the centered wrapper associated with $f$.
\end{prop}
\begin{figure}[h]
\begin{centering}
\includegraphics[width=7cm]{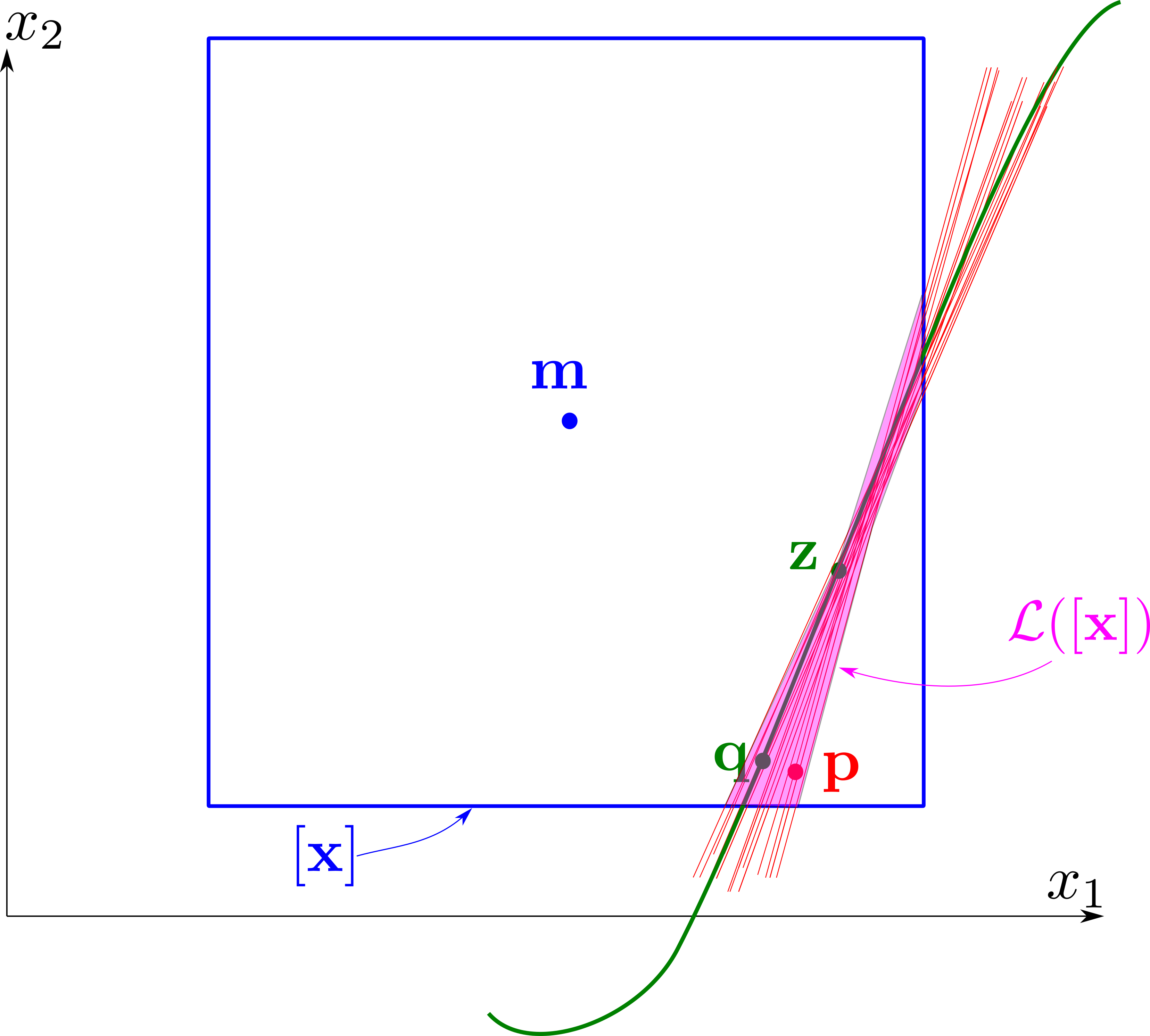}
\par\end{centering}
\caption{The set $\mathcal{L}([\mathbf{x}])$ (red) with a bowtie shape is
close to the curve $\mathbb{X}$ (green)}
\label{fig:interval:lin:approx}
\end{figure}

\textbf{Proof}. Consider the sequence $[\mathbf{x}](k)\subset\mathbb{R}^{n}$
converging to $\mathbf{z}$. We assume that $[\mathbf{x}](k)$, or
$[\mathbf{x}]$ for short, is narrow, \emph{i.e.}, $w([\mathbf{x}])=O(\varepsilon)$.
If $\mathbf{p}\in\mathcal{L}([\mathbf{x}])$ (see Figure \ref{fig:interval:lin:approx})
then,~for some $\mathbf{a}\in[\mathbf{a}]=[\frac{df}{d\mathbf{x}}]([\mathbf{x}])$,
we have
\begin{equation}
f(\mathbf{m})+\mathbf{a}\cdot(\mathbf{p}-\mathbf{m})=0
\end{equation}
where $\mathbf{m}=\text{center}(\text{[\ensuremath{\mathbf{x}}]})$.
From Corollary \ref{cor:moveplane}, taking $d\mathbf{x}=\mathbf{p}-\mathbf{m}=O(\varepsilon)$
and since $w([\mathbf{a}])=O(\varepsilon)$, we get that the distance
between a point in $\mathcal{L}([\mathbf{x}])$ and the set $\mathbb{X}$
is an $o(\varepsilon)$. We get that
\begin{equation}
h(\mathcal{L}([\mathbf{x}](k)),\mathbb{X})=o(w([\mathbf{x}](k)))
\end{equation}
i.e.,
\begin{equation}
\lim_{k\rightarrow\infty}\frac{h(\mathcal{L}([\mathbf{x}](k)),\mathbb{X})}{w([\mathbf{x}](k))}=0.
\end{equation}
Thus the wrapper $\mathcal{L}$ is of order 1 at $\mathbf{z}$.$\blacksquare$
\begin{cor}
The contractor for $f(\mathbf{x})=0$ defined by
\begin{equation}
\begin{array}{ccl}
[x_{i}] & = & [x_{i}]\cap\left(m_{i}-f(\mathbf{m})-\sum_{j\neq i}[a_{j}]\cdot([x_{j}]-m_{j})\right)\\{}
[a_{j}] & = & [\frac{\partial f}{\partial x_{j}}]([\mathbf{x}])
\end{array}
\end{equation}
is asymptotically minimal.
\end{cor}
\textbf{Proof}. Define $\mathcal{L}([\mathbf{x}])$ as in (\ref{eq:Lxx}).
From Proposition \ref{prop:Cmini:fromW}, $\mathcal{L}\in\text{Wrap}(\mathbb{X},\mathbf{z}).$
The contractor $\mathcal{C}([\mathbf{x}])=\llbracket\mathcal{L}([\mathbf{x}])\rrbracket$
is an asymptotically minimal contractor. Now the set $\mathcal{L}([\mathbf{x}])$
can be defined by the following constraints
\begin{equation}
\left\{ \begin{array}{c}
\exists\mathbf{z}\in[\mathbf{x}]\\
f(\mathbf{m})+\mathbf{a}\cdot(\mathbf{x}-\mathbf{m})=0\\
\mathbf{a}=\frac{\partial f}{\partial\mathbf{x}}(\mathbf{z})\\
\mathbf{m}=\text{center}([\mathbf{x}])
\end{array}\right.
\end{equation}
Since $\mathbf{x}$ occurs only once in the constraint $f(\mathbf{m})+\mathbf{a}\cdot(\mathbf{x}-\mathbf{m})=0$,
an interval forward-backward propagation provides us the minimal contraction
\cite{Montanari91}, \emph{i.e.}, it returns the box $\llbracket\mathcal{L}([\mathbf{x}])\rrbracket$.$\blacksquare$

\subsection{Vector case}
\begin{prop}
\label{prop:Lxx:vect}Consider the equation $\mathbf{f}(\mathbf{x})=\mathbf{0}$,
where \textbf{$\mathbf{f}:\mathbb{R}^{n}\mapsto\mathbb{R}^{p}$} is
differentiable. The solution set is
\begin{equation}
\begin{array}{ccl}
\mathbb{X} & = & \{\mathbf{x}\in\mathbb{R}^{n}\,|\,\mathbf{f}(\mathbf{x})=\mathbf{0}\}.\end{array}
\end{equation}
Consider a point $\mathbf{z}$ such that $\mathbf{f}(\mathbf{z})=\mathbf{0}$
and a nested sequence $[\mathbf{x}](k)$ converging to $\mathbf{z}$.
Consider the wrappers $\mathcal{L}_{i}:\mathbb{IR}^{n}\mapsto\mathcal{P}(\mathbb{R}^{n})$
of order 1 for $f_{i}(\mathbf{x})=0$ defined by
\begin{equation}
\begin{array}{ccc}
\mathcal{L}_{i}([\mathbf{x}]) & =\{ & \mathbf{x}\in[\mathbf{x}]\,|\,\exists\mathbf{a}\in[\frac{df_{i}}{d\mathbf{x}}]([\mathbf{x}]),\\
 &  & f_{i}(\mathbf{m})+\mathbf{a}_{i}\cdot(\mathbf{x}-\mathbf{m})=0,\\
 &  & \mathbf{m}=\text{center}([\mathbf{x}])\}.
\end{array}\label{eq:Lxx:vect}
\end{equation}
The operator $\bigcap_{i}\mathcal{L}_{i}$, belongs to $\text{Wrap}(\mathbb{X},\mathbf{z})$. 
\end{prop}
\textbf{Proof}. It is a direct consequence of Proposition \ref{prop:n:wrappers}.
$\blacksquare$

To compute $\bigcap_{i}\mathcal{L}_{i}$, the method proposed for
the scalar case is not valid anymore. An interval linear method could
be used \cite{DBLP:Neumaier:Shcherbina} \cite{Araya12}. An other
possibility is to use a preconditioning method based on the Gauss-Jordan
decomposition, which will be minimal in many cases, such as the test-case
that will be treated in Section \ref{sec:testcase}. 

\subsection{Preconditioning}

Consider the equation $\mathbf{f}(\mathbf{x})=\mathbf{0}$, where
\textbf{$\mathbf{f}:\mathbb{R}^{n}\mapsto\mathbb{R}^{p}$} is differentiable.
Intersecting sets $\mathcal{L}_{i}([\mathbf{x}])$ as suggested by
Proposition \ref{prop:Lxx:vect} requires the resolution of interval
linear equations. This operation is costly and should be avoided if
it has to be repeated a large number of times. Instead of this, we
prefer to use a specific preconditioning method.

\begin{figure}[h]
\begin{centering}
\centering\includegraphics[width=5cm]{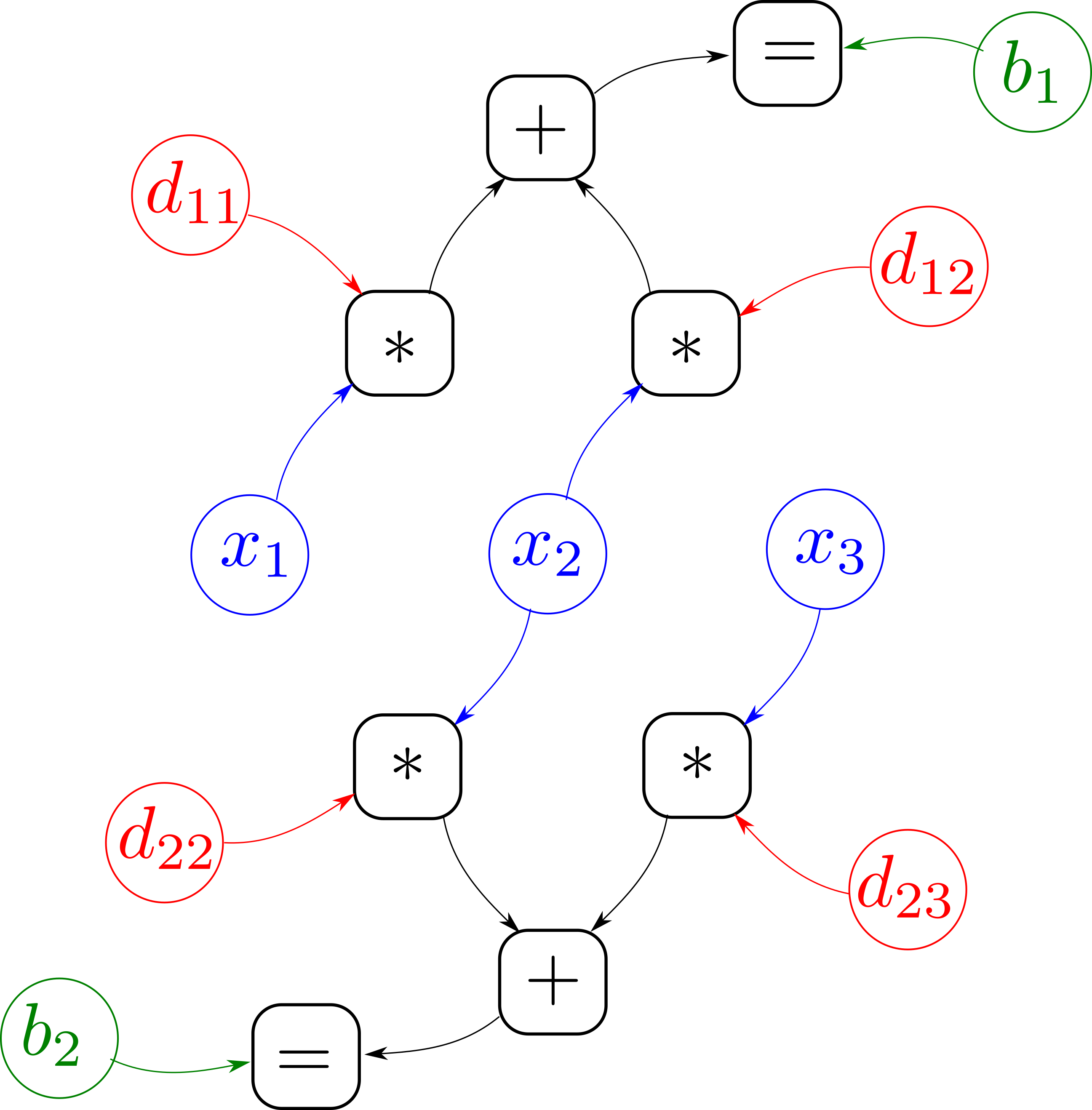}
\par\end{centering}
\caption{The constraint network has no cycle (it is a tree). Thus the interval
propagation is minimal}
\label{fig:lu.png}
\end{figure}

To understand the principle of the preconditioning, consider the following
interval linear system 
\begin{equation}
\left(\begin{array}{ccc}
d_{11} & d_{12} & 0\\
0 & d_{22} & d_{23}
\end{array}\right)\left(\begin{array}{c}
x_{1}\\
x_{2}\\
x_{3}
\end{array}\right)=\left(\begin{array}{c}
b_{1}\\
b_{2}
\end{array}\right)\label{eq:lu23}
\end{equation}
where
\begin{equation}
d_{ij}\in[d_{ij}],x_{j}\in[x_{j}],b_{i}\in[b_{i}]
\end{equation}
The optimal contraction can be obtained by a simple interval propagation.
This is due to the fact that the corresponding constraint network
as no cycle \cite{Montanari91}, as illustrated by Figure \ref{fig:lu.png}. 

Note that no cycle would have been obtained with the following linear
system:
\begin{equation}
\left(\begin{array}{cccc}
d_{11} & d_{12} & 0 & 0\\
0 & d_{22} & d_{23} & 0\\
0 & 0 & d_{33} & d_{34}
\end{array}\right)\left(\begin{array}{c}
x_{1}\\
x_{2}\\
x_{3}\\
x_{4}
\end{array}\right)=\left(\begin{array}{c}
b_{1}\\
b_{2}\\
b_{3}
\end{array}\right)\label{eq:lu34}
\end{equation}

A matrix $\mathbf{D}$ such that the system $\mathbf{D}\cdot\mathbf{x}=\mathbf{b}$
has no cycle can be called a \emph{tree matrix}. 

Both systems (\ref{eq:lu23}) and (\ref{eq:lu34}), for which the
matrix $\mathbf{D}$ is a \emph{band }matrix \cite{Atkinson89}, could
be obtained from a Gauss Jordan transformation of a linear systems
\cite{Leon09}. For instance, if we have a system of the form $\mathbf{Ax}=\mathbf{c}$
where $\mathbf{A}$ is of dimension $3\times4$ with full rank, there
exists a matrix $\mathbf{Q}$ of dimension $3\times3$ such that
\begin{equation}
\mathbf{Ax}=\mathbf{c}\Leftrightarrow\mathbf{Q}\cdot\mathbf{A\cdot}\mathbf{x}=\mathbf{Q}\cdot\mathbf{c}
\end{equation}
where $\mathbf{D}=\mathbf{Q}\cdot\mathbf{A}$ has the form given by
(\ref{eq:lu34}). 
\begin{prop}
\label{prop:treematrix:cform}Consider a set $\mathbb{X}=\left\{ \mathbf{x}\in\mathbb{R}^{n}|\mathbf{f}(\mathbf{x})=\mathbf{0}\right\} $.
Take a narrow box $[\mathbf{x}]$ with center $\mathbf{m}$. Assume
that $\frac{d\mathbf{f}}{d\mathbf{x}}(\mathbf{m})$ is a tree matrix.
An interval propagation on the system
\begin{equation}
\begin{array}{c}
\mathbf{f}(\mathbf{m})+\mathbf{A}\cdot(\mathbf{x}-\mathbf{m})=\mathbf{0}\\
\mathbf{A}\in[\frac{d\mathbf{f}}{d\mathbf{x}}]([\mathbf{x}])\\
\mathbf{x}\in[\mathbf{x}]
\end{array}
\end{equation}
corresponds to an asymptotically minimal contractor for $\mathbb{X}.$
\end{prop}
\textbf{Proof.} The interval matrix $[\mathbf{A}]=[\frac{d\mathbf{f}}{d\mathbf{x}}]([\mathbf{x}])$
is such that $w([\mathbf{A}])=O(\varepsilon)$, where $\varepsilon=w([\mathbf{x}])$.
Due to the fact that the contractor $\mathcal{C}$ resulting from
the interval propagation is minimal for $\mathbf{A=}\frac{d\mathbf{f}}{d\mathbf{x}}(\mathbf{m})$,
and taking into account Proposition \ref{prop:sensitivity:lin}, we
get that the contractor obtained by an elementary interval propagation
is asymptotically minimal.$\blacksquare$ 
\begin{cor}
Consider a set $\mathbb{X}=\left\{ \mathbf{x}\in\mathbb{R}^{n}|\mathbf{f}(\mathbf{x})=\mathbf{0}\right\} $.
Take a narrow box $[\mathbf{x}]$ with center $\mathbf{m}$. Define
$\mathbf{Q}$ such that $\mathbf{Q}\cdot\frac{d\mathbf{f}}{d\mathbf{x}}(\mathbf{m})$
is a tree matrix. An interval propagation on the system
\begin{equation}
\begin{array}{c}
\mathbf{Q}\cdot\mathbf{f}(\mathbf{m})+\mathbf{Q}\cdot\mathbf{A}\cdot(\mathbf{x}-\mathbf{m})=\mathbf{0}\\
\mathbf{A}\in[\frac{d\mathbf{f}}{d\mathbf{x}}]([\mathbf{x}])\\
\mathbf{x}\in[\mathbf{x}]
\end{array}
\end{equation}
corresponds to an asymptotically minimal contractor for $\mathbb{X}.$ 
\end{cor}
\textbf{Proof.} It suffices to apply Proposition \ref{prop:treematrix:cform}
with $\mathbf{g}(\mathbf{x})=\mathbf{Q}\cdot\mathbf{f}(\mathbf{x})$.$\blacksquare$

\subsection{Algorithm\label{subsec:algorithm}}

Consider the system $\mathbf{f}(\mathbf{x})=\mathbf{0}$ and take
a box $[\mathbf{x}]$. The following algorithm corresponds to a centered
contractor.

\begin{tabular}{|cl|}
\hline 
Input:  & \textbf{f},$[\mathbf{x}]$\tabularnewline
\hline 
\hline 
1 & $\mathbf{m}=\text{center}([\mathbf{x}])$ \tabularnewline
2 & Compute the Gauss-Jordan matrix $\mathbf{Q}$ for $\frac{d\mathbf{f}}{d\mathbf{x}}(\mathbf{m})$ \tabularnewline
3 & Define $\mathbf{g}(\mathbf{x})=\mathbf{Q}\cdot\mathbf{f}(\mathbf{x})$\tabularnewline
4 & For $i\in\{1,\dots,p\}$\tabularnewline
5 & \quad{}For $j\in\{1,\dots,n\}$\tabularnewline
6 & \quad{}\quad{}$[\mathbf{a}]=$$[\frac{\partial g_{i}}{\partial\mathbf{x}}]([\mathbf{x}])$\tabularnewline
7 & \quad{}\quad{}$[s]={\displaystyle \sum_{k\neq j}}[a_{k}]\cdot([x_{k}]-m_{k})$\tabularnewline
8 & \quad{}\quad{}$[x_{j}]=[x_{j}]\cap\left(-g_{i}(\mathbf{m})-[s]\right)$\tabularnewline
9 & Return $[\mathbf{x}]$\tabularnewline
\hline 
\end{tabular}
\begin{itemize}
\item Step 1 takes the center $\mathbf{m}$ of $[\mathbf{x}]$ in order
to form a linear approximation for $\mathbf{f}$ in $[\mathbf{x}]$:
\begin{equation}
\mathbf{f}(\mathbf{x})=\mathbf{f}(\mathbf{m})+\frac{d\mathbf{f}}{d\mathbf{x}}(\mathbf{m})\cdot(\mathbf{x}-\mathbf{m}).
\end{equation}
\item Step 2 returns an invertible $m\times m$ matrix $\mathbf{Q}$ such
that $\mathbf{A}=\mathbf{Q}\cdot\frac{d\mathbf{f}}{d\mathbf{x}}(\mathbf{m})$
is a band matrix. The matrix $\mathbf{Q}$ is chosen by a Gauss-Jordan
algorithm. The new system to be solved is now
\begin{equation}
\mathbf{Q}\cdot\mathbf{f}(\mathbf{x})=\mathbf{0}.
\end{equation}
\item Step 3 defines $\mathbf{g}(\mathbf{x})=\mathbf{Q}\cdot\mathbf{f}(\mathbf{x}-\mathbf{m}).$
We need to solve $\mathbf{g}(\mathbf{x})=\mathbf{0}$ in the box $[\mathbf{x}]-\mathbf{m}$.
The main difference compared to the previous system $\mathbf{f}(\mathbf{x})=\mathbf{0}$
is that its linear approximation
\begin{equation}
\mathbf{g}(\mathbf{x})=\mathbf{g}(\mathbf{m})+\mathbf{A}\cdot(\mathbf{x}-\mathbf{m}).
\end{equation}
is such that $\mathbf{A}$ is a band matrix.
\item Step 4-9 define the set of constraints
\begin{equation}
\left\{ \begin{array}{c}
\mathbf{0}=\mathbf{g}(\mathbf{m})+\mathbf{A}\cdot(\mathbf{x}-\mathbf{m})\\
\mathbf{A}\in[\frac{d\mathbf{g}}{d\mathbf{x}}]([\mathbf{x}])\\
\mathbf{x}\in[\mathbf{x}]
\end{array}\right.
\end{equation}
 and performs an interval propagation. Due to the fact that the system
has no cycle (at first order) then the propagation is asymptotically
minimal. 
\end{itemize}

\section{Test case\label{sec:testcase}}

Interval methods have been shown to be very powerful for the stability
analysis of linear systems \cite{Malan96}. We have chosen to consider
the linear time-delay system \cite{turkulov:hal-03184836} given by
\begin{equation}
\ddot{x}+2\dot{x}(t-p_{1})+x(t-p_{2})=0
\end{equation}
but other types of linear systems \cite{DBLP:VictorMMCA22} with fractional
orders could be considered as well. Its characteristic function is

\begin{equation}
\theta(\mathbf{p},s)=s^{2}+2se^{-sp_{1}}+e^{-sp_{2}}.
\end{equation}
For a given $\mathbf{p}$, the location of the roots for $\theta(\mathbf{p},s)$
provides an information concerning the stability of the system. For
instance, if all roots are on the half left of the complex plane,
then the system is stable. The stability changes when one root crosses
the imaginary line. This is the reason why we are interested in characterizing
the set
\begin{equation}
\mathcal{P}=\{\mathbf{p}\,|\,\exists\omega>0,\,\theta(\mathbf{p},j\omega)=0\}.
\end{equation}

Now
\begin{equation}
\begin{array}{cl}
 & \theta(p_{1},p_{2},j\omega)\\
= & -\omega^{2}+2j\omega e^{-j\omega p_{1}}+e^{-j\omega p_{2}}\\
= & -\omega^{2}+2j\omega(\cos(\omega p_{1})-j\sin(\omega p_{1}))\\
 & +\cos(\omega p_{2})-j\sin(-\omega p_{2})\\
= & -\omega^{2}+2\omega\sin(\omega p_{1})+\cos(\omega p_{2})\\
 & +j\cdot\left(2\omega\cos(\omega p_{1})-\sin(\omega p_{2})\right)
\end{array}
\end{equation}
We have
\begin{equation}
\begin{array}{cl}
 & \theta(p_{1},p_{2},j\omega)=0\\
\Leftrightarrow & \underset{\mathbf{f}(p_{1},p_{2},\omega)}{\underbrace{\left(\begin{array}{c}
-\omega^{2}+2\omega\sin(\omega p_{1})+\cos(\omega p_{2})\\
2\omega\cos(\omega p_{1})-\sin(\omega p_{2})
\end{array}\right)}}=0
\end{array}
\end{equation}

\begin{figure}[h]
\begin{centering}
\centering\includegraphics[width=7cm]{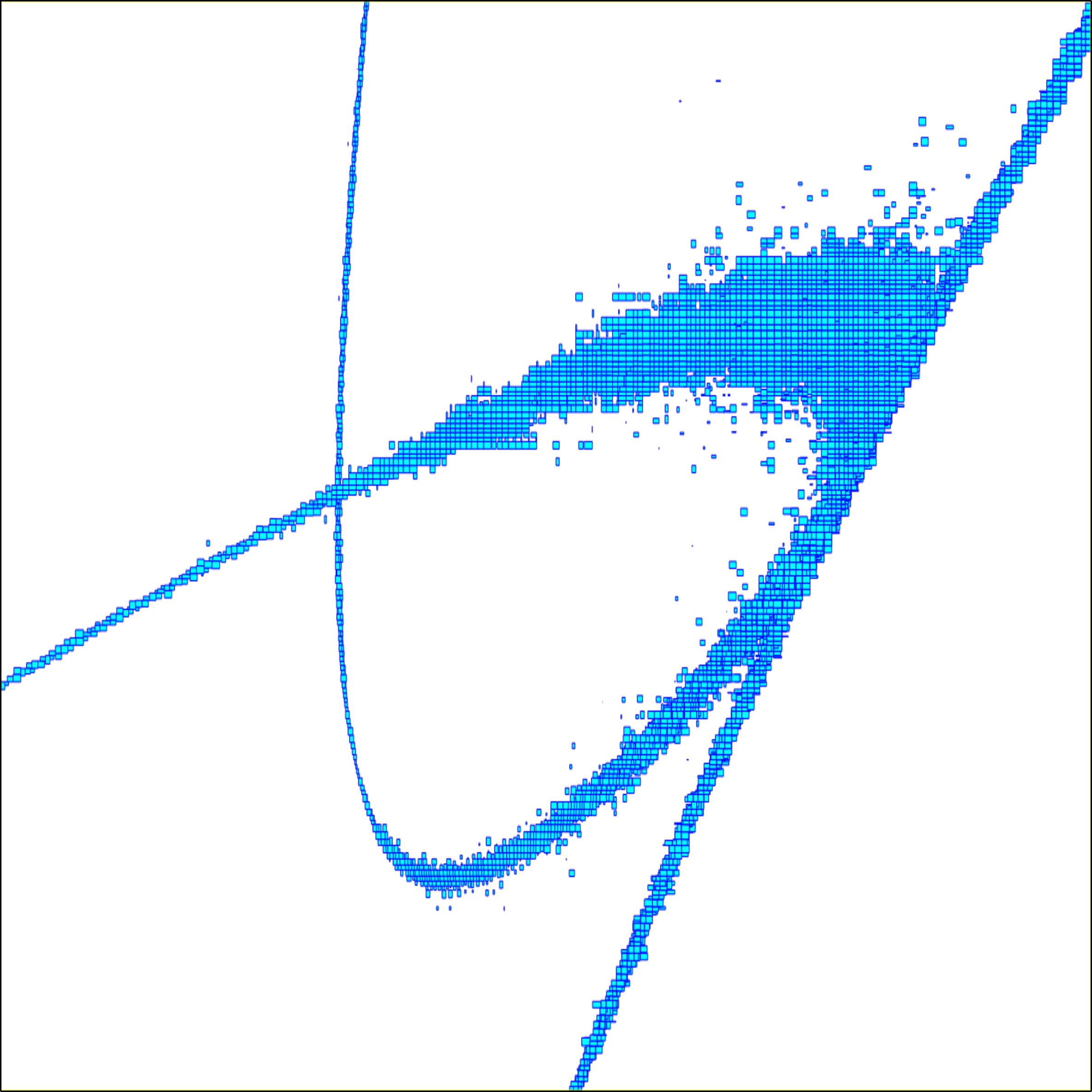}
\par\end{centering}
\caption{Approximation of the solution set with a state of the art contractor.
The frame box for $(p_{1},p_{2})$ is $[0,2.5]\times[2,4]$}
\label{fig:paving0}
\end{figure}

Take $[p_{1}]=[0,2.5]$, $[p_{2}]=[1,4]$,$[\omega]=[0,10]$ and let
us characterize the set $\mathcal{P}$ using the centered contractor.
Using a branch and prune algorithm with a accuracy of $\varepsilon=2^{-8}$
with an HC4 algorithm \cite{Ceberio01}\cite{BenhamouICLP99} (the
state of the art), we get the paving of Figure \ref{fig:paving0}
in 4 sec. The number of boxes of the approximation is $43173.$ Similar
results where obtained were obtained on the same example in \cite{malti:hal-03646956}.

\begin{figure}[h]
\begin{centering}
\centering\includegraphics[width=7cm]{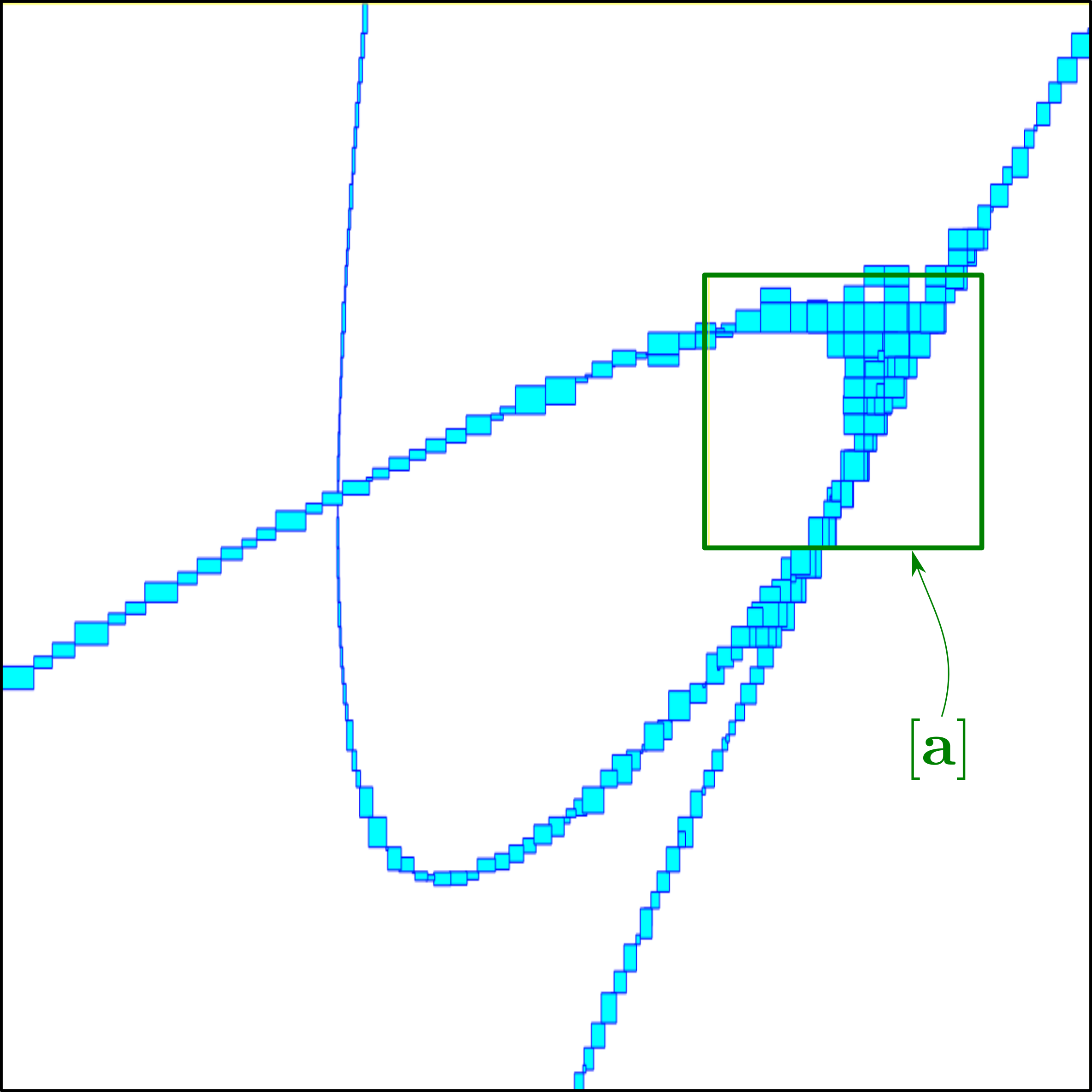}
\par\end{centering}
\caption{Paving obtained with the centered contractor. The frame box for $(p_{1},p_{2})$
is $[0,2.5]\times[2,4]$}
\label{fig:paving1}
\end{figure}

With an accuracy of $\varepsilon=2^{-4}$ with the centered contractor
given in Section \ref{subsec:algorithm}, we get the paving of Figure
\ref{fig:paving1} in 1.2 sec. The number of boxes of the approximation
is 282 (instead of 43173), for a more accurate approximation. 

\begin{figure}[h]
\begin{centering}
\centering\includegraphics[width=7cm]{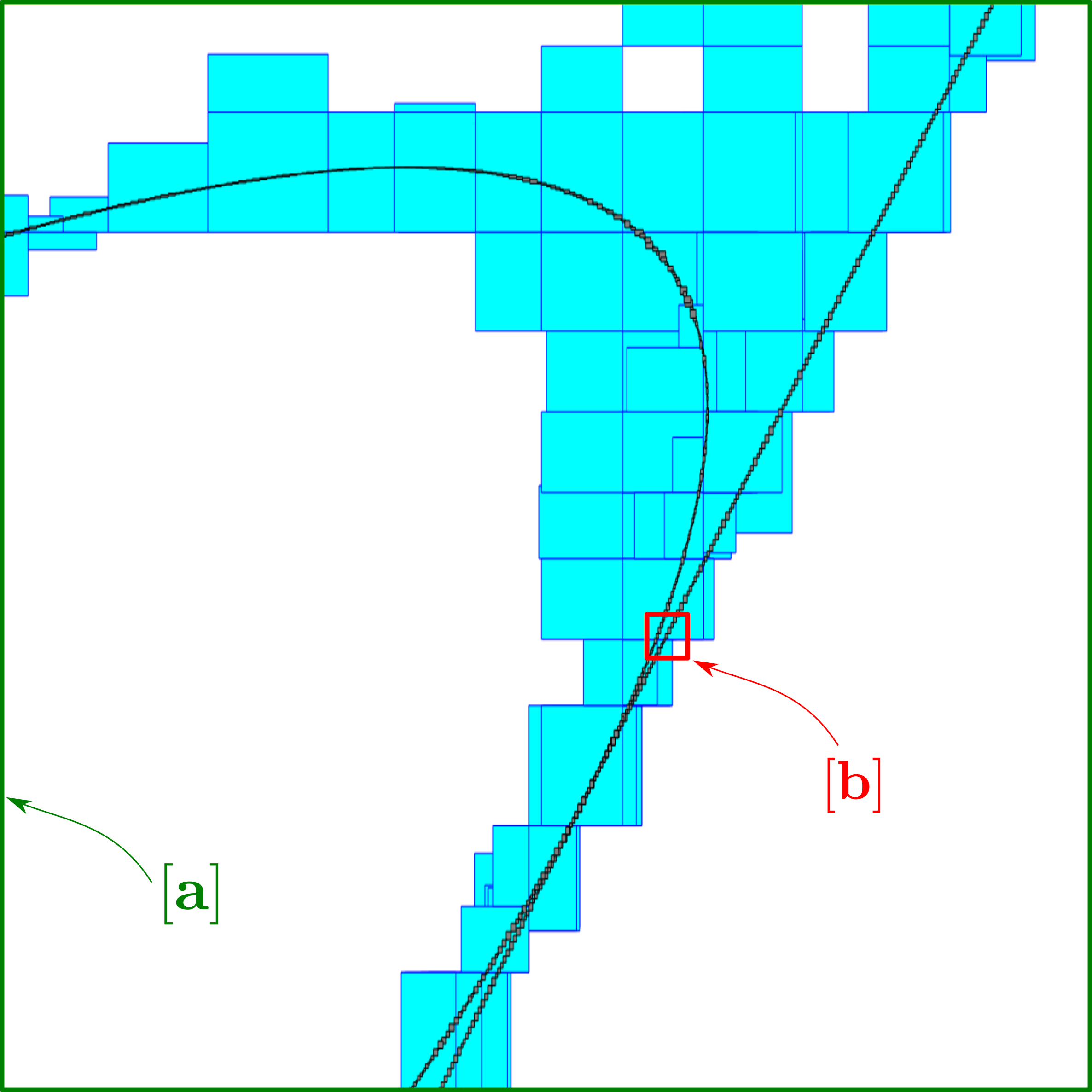}
\par\end{centering}
\caption{Pavings obtained with the centered contractor in the box $[\mathbf{a}]=[1.3,1.8]\times[3.0,3.5]$;
Blue: $\varepsilon=2^{-4}$ ; Thin: $\varepsilon=2^{-8}$ }
\label{fig:paving2}
\end{figure}

With a accuracy of $\varepsilon=2^{-8}$ with the centered contractor,
we get the thin curve represented on Figure \ref{fig:paving2}. This
curve is made with the small boxes generated by the paver, which shows
the quality of the approximation. The big blue boxes are those already
painted in the green box $[\mathbf{a}]$ of Figure \ref{fig:paving1}. 

\begin{figure}[h]
\begin{centering}
\centering\includegraphics[width=7cm]{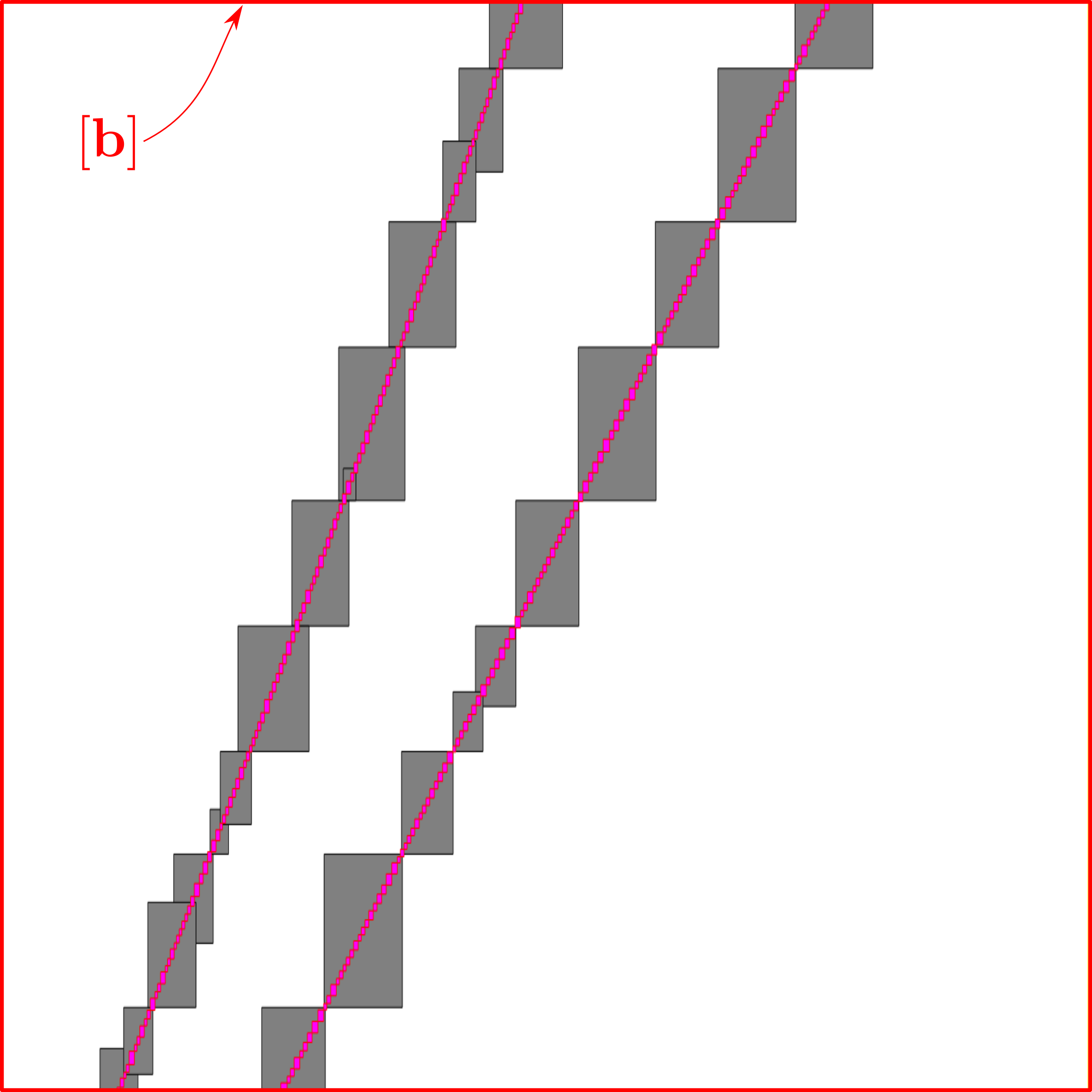}
\par\end{centering}
\caption{Approximation of the solution set in $[\mathbf{b}]=[1.595,1.615]\times[3.2,3.22]$;
Gray: $\varepsilon=2^{-8}$ ; Magenta: $\varepsilon=2^{-12}$}
\label{fig:paving3}
\end{figure}

With a accuracy of $\varepsilon=2^{-12}$ with the centered contractor,
we get the magenta curve of Figure \ref{fig:paving3}. The big gray
boxes are those already painted in the red box $[\mathbf{b}]$ of
Figure \ref{fig:paving2}. The fact that, for a small $\varepsilon$
, the boxes of the approximation only overlap on their corners illustrates
the minimality of the contractor. 

The computing time to get the three Figures \ref{fig:paving1}, \ref{fig:paving2}
and \ref{fig:paving3} is less than 10 sec. Our results are much more
accurate than those obtained in Section 6 of \cite{malti:hal-03646956}. 

The code and an illustrating video are given at \\
\texttt{www.ensta-bretagne.fr/jaulin/centered.html}

\section{Conclusion\label{sec:Conclusion}}

In this paper, we have proposed a contractor which is asymptotically
minimal for the approximation of a curve defined by nonlinear equations.
The resulting \emph{centered }contractor is based on the centered
form which suppresses the pessimism when the boxes are narrow and
when we have a single equation. When we combine several equations,
a preconditioning method has been proposed in order to linearize the
problem into a system where a tree matrix in involved. The preconditioning
has been implemented using a Gauss Jordan band diagonalization method.
On an example, we have shown that our centered contractor was able
to outperform the state of the art contractor based on a forward-backward
propagation.

Other approaches, such as the generalized interval arithmetic \cite{Hansen75},
the affine arithmetic \cite{Figueiredo04} allows to get first order
approximation of the constraints. As for our paper, these arithmetics
can obviously model the affine dependencies between quantities with
an error that shrinks quadratically with the size of the input intervals.
Now, this linear approximation is only valid when we have a single
constraint and can thus not be used to build asymptotically minimal
contractors without some improvements. Our approach 
\begin{itemize}
\item does not require the implementation of a new arithmetic; it only uses
the standard interval arithmetic 
\item generates a contractor that can be combined with other existing contractors
enforcing the efficiency of the resolution.
\end{itemize}


\bibliographystyle{plain}

\end{document}